\documentclass[reqno]{amsart}
\usepackage{hyperref}
\usepackage{amsmath}
\usepackage{amssymb}
\usepackage{amscd}
\usepackage{graphics}
\usepackage{latexsym}
\usepackage{stmaryrd}

\theoremstyle{plain}
\newtheorem{theorem}{Theorem}[section]
\newtheorem{thm}[theorem]{Theorem}
\newtheorem{cor}[theorem]{Corollary}
\newtheorem{lem}[theorem]{Lemma}
\newtheorem{prop}[theorem]{Proposition}

\theoremstyle{definition}
\newtheorem{defn}[theorem]{Definition}
\newtheorem{prob}[theorem]{Problem}

\newtheorem{rmk}[theorem]{Remark}

\newtheorem{notat}[theorem]{Notation}

\newtheorem{hyp}[theorem]{Hypothesis}

\theoremstyle{remark}

\newcommand{\ZZ}{\mathbb{Z}}

\newcommand{\RR}{\mathbb{R}}

\newcommand{\PP}{\mathbb{P}}

\newcommand{\mc}{\mathcal}

\newcommand{\OO}{\mc{O}}

\newcommand{\marpar}[1]{}

\newcommand{\lt}{\left}
\newcommand{\rt}{\right}

\newsavebox{\sembox}
\newlength{\semwidth}
\newlength{\boxwidth}

\newsavebox{\semrbox}
\newlength{\semrwidth}
\newlength{\boxrwidth}

\title{A note on Fano manifolds whose second Chern character is positive} 
\author[de Jong]{A. J. de Jong}
\author[Starr]{Jason Michael Starr}
\begin{document}


\begin{abstract}
This note outlines some first steps in the classification of Fano
manifolds for which $c_1^2-2c_2$ is positive or nef.
\end{abstract}


\maketitle


\section{Introduction} \label{sec-intro}

\noindent
This note about Fano manifolds $X$ for which
$(c_1^2-2c_2)(T_X)$ is positive, lists what few examples are known, as well as
giving many non-examples.  
Presumably there are many more examples.
They do not seem easy to find.

\begin{notat} \label{notat-N}
Let $X$ be a projective variety over an algebraically closed
field.  For every integer $k\geq 0$, denote by $N_k(X)$ the
finitely-generated free Abelian group of $k$-cycles modulo numerical
equivalence, and denote by $N^k(X)$ the $k^\text{th}$ graded piece of the
quotient algebra $A^*(X)/\text{Num}^*(X)$, cf. ~\cite[Example
  19.3.9]{F}.  
For every $\ZZ$-module $B$, denote $N_k(X)_B :=
N_k(X)\otimes B$, resp. $N^k(X)_B := N^k(X)\otimes B$.  
Denote by $NE_k(X) \subset N_k(X)$ the semigroup generated by
nonzero, effective $k$-cycles.  
For $B$ a subring of $\RR$, denote by
$NE_k(X)_B$ the $B_{>0}$-semigroup in $N_k(X)_B$ generated by $NE_k(X)$.
\end{notat}

\begin{defn} \label{defn-N}
A class in $N^k(X)_\RR$ is \emph{nef} if it pairs nonnegatively
with every element in $\overline{NE}_k(X)$.  The corresponding cone is
denote $\text{Nef}^k(X)$.  
A class is
\emph{weakly positive} if it pairs positively with every element in
$NE_k(X)$.  The corresponding cone is denoted $\text{WPos}^k(X)$.
A class is \emph{positive} if it is contained
in the interior of $\text{Nef}^k(X)$; the interior of
$\text{Nef}^k(X)$ is denoted $\text{Pos}^k(X)$.  The \emph{ample cone}
is the $\RR_{>0}$-semigroup generated by the
image of the cup-product map, $(\text{Pos}^1(X))^k \rightarrow
N^k(X)_\RR$.  It is denoted $\text{Ample}^k(X)$, and its elements are
\emph{ample} classes.
\end{defn}

\begin{rmk} \label{rmk-N}
There are obvious inclusions, 
$$
\text{Ample}^k(X) \subset \text{Pos}^k(X) \subset \text{WPos}^k(X)
\subset \text{Nef}^k(X).
$$
For $k=1$, $\text{Ample}^1(X) = \text{Pos}^1(X)$ by definition.
Moreover, by Kleiman's criterion, this is the $\RR_{>0}$-semigroup
generated by first Chern classes of ample invertible sheaves.  For
$k>1$, it can happen that $\text{Ample}^k(X) \neq \text{Pos}^k(X)$;
for instance, because $(N^1(X))^{\otimes k} \rightarrow N^k(X)$ is not
surjective.  There are also examples where $\text{Pos}^k(X) \neq
\text{WPos}^k(X)$ and $\text{WPos}^k(X) \neq \text{Nef}^k(X)$.
\end{rmk}

\begin{prob} \label{prob-main}
Find smooth, connected, projective varieties $X$ such that
$\text{ch}_1(T_X)=c_1(T_X)$ is ample and 
$\text{ch}_2(T_X)= 1/2(c_1^2-2c_2)(T_X)$ is ample, resp. positive, weakly
positive, nef.  More generally, allow $X$ to be a
smooth, connected, proper Deligne-Mumford stack whose coarse moduli
space is projective.
\end{prob}

\section{Positive Examples} \label{sec-pexamples}

\medskip\noindent
Following are examples of Fano manifolds with $\text{ch}_2(T_X)$ ample
or positive.

\medskip\noindent
\textbf{1.}
The simplest example is $\PP^n$ for $n\geq 2$.  Denote by $h\in
N^1(\PP^n)$ the first Chern class of $\OO_{\PP^n}(1)$.  Using the
Euler sequence,
$$
\begin{CD}
0 @>>> \OO_{\PP^n} @>>> \OO_{\PP^n}(1)^{\oplus (n+1)} @>>> T_{\PP^n}
@>>> 0,
\end{CD}
$$
the Chern character of $T_{\PP^n}$ is $(n+1)e^h - 1$.  In particular,
$\text{ch}_k(T_X) = (n+1)h^k/k!$ for every $k=1,\dots,n$.  
So $\text{ch}_k(T_X)$ is ample for $k=1,\dots,n$.

\medskip\noindent
\textbf{2.}
Weighted projective spaces are also examples.  
The weighted projective space $\PP = \PP(a_0,\dots,a_n)$ is the coarse
moduli space of a smooth Deligne-Mumford stack $X$,
and
$\text{ch}_k(T_X) = (n+1)h^k/k!$ where
$h$ is the first Chern class of the invertible sheaf $\OO_{X}(1)$ on
the stack.  Some positive multiple of $h$ is the pullback of an ample
class from the coarse moduli space, thus $h$ is an ample class.  

\medskip\noindent
\textbf{3.}
Let $Y$ be a smooth complete intersection of divisors $D_1,\dots,D_r$
in $\PP$ of respective degrees $d_1,\dots,d_r$.  Using the exact
sequences,
$$
\begin{CD}
0 @>>> T_Y @>>> T_{\PP}|_Y @>>> \oplus_{i=1}^r \OO_{\PP}(d_i)|_Y @>>>
0,
\end{CD}
$$
the Chern character of $T_Y$ is $(n+1)e^h - 1 - \sum_{i=1}^r e^{d_i
  h}$.  Thus $\text{ch}_k(T_Y) = 1/k!(n+1 - (d_1^k + \dots +
d_r^k))h^k$ for $k=1,\dots,n-r$.  In particular, if $d_1^2 + \dots +
d_r^k < n+1$ then
$\text{ch}_1(T_Y)$ and $\text{ch}_2(T_Y)$ are both ample.

\medskip\noindent
\textbf{4.}
For every integer $k\geq 1$, the Grassmannians
$G=\text{Grass}(k,2k)$ and $G=\text{Grass}(k,2k+1)$ have
$\text{ch}_1(T_G)$ is
ample and $\text{ch}_2(T_G)$ is positive.  If $k>1$, then 
$\text{ch}_2(T_G)$ is not in the $1$-dimensional subspace
spanned by $\sigma_1^2$.  Therefore $\text{ch}_2(T_G)$ is positive,
but not ample.

\section{Nef examples} \label{sec-nefexamples}

\medskip\noindent
Given a Fano manifold, there are a few
methods of constructing a new Fano manifold $Y$ with $\text{ch}_2(T_Y)$
nef.  Typically even if $\text{ch}_2(T_X)$ is positive,
$\text{ch}_2(T_Y)$ is not weakly positive.

\medskip\noindent
\textbf{1.}
Let $X$ be a Fano manifold with $\text{ch}_2(T_X)$ nef.  Let $Y$ be a
smooth divisor in $X$.  If $c_1(T_X) - [Y]$ is ample, and
$\text{ch}_2(T_X) - [Y]^2/2$ is nef, then $Y$ is a Fano manifold and
$\text{ch}_2(T_Y)$ is nef.  This is essentially the same as Example 3
in Section~\ref{sec-pexamples}.

\medskip\noindent
\textbf{2.}
Let $X$ be a Fano manifold and let $L$ be a nef line bundle such that
$c_1(T_X) - c_1(L)$ is ample and $\text{ch}_2(T_X) + c_1(L)^2/2$ is
nef.  Then the projective bundle $\PP E = \PP(L^\vee \oplus \OO_X)$ is
a Fano manifold and $\text{ch}_2(T_{\PP E})$ equals
$\pi^*(\text{ch}_2(T_X) + c_1(L)^2/2)$.  This is nef, but not weakly
positive; its restriction to $\pi^{-1}(C)$ is zero for every curve
$C\subset X$.
Note $\text{ch}_2(T_X)$ need not be nef, e.g., for integers $(n,d,a)$
satisfying
$1\leq d\leq
\lfloor (n^2+n+1)/2n \rfloor$ and $\lceil \sqrt{\max(0,d^2-n-1)} \rceil
\leq a\leq n-d$,  
for every smooth
degree $d$ hypersurface $X\subset \PP^n$, the projective bundle
$\PP(\OO_X(-a)\oplus \OO_X)$ is a Fano manifold with
$\text{ch}_2(T_{\PP E})$ nef.  

\medskip\noindent
\textbf{3.}  
Let $X$ and $Y$ be Fano manifolds such that $\text{ch}_2(T_X)$ and
$\text{ch}_2(T_Y)$ are nef.  The product $X\times Y$ is Fano and
$\text{ch}_2(T_{X\times Y}) =
\pi_X^*\text{ch}_2(T_X)+\pi_Y^*\text{ch}_2(T_Y)$, which is nef.
For rational curves $C_X\subset X$ and $C_Y\subset Y$, the pairing of
$\text{ch}_2(T_{X\times Y})$ with $C_X\times C_Y$ is zero, thus
$\text{ch}_2(T_{X\times Y})$ is not weakly positive.

\section{Projective bundles} \label{sec-pbundles}

\medskip\noindent
One way to produce new examples of Fano manifolds is to form the
projective bundle of a vector bundle of ``low degree'' over a given
Fano manifold.  

\begin{lem} \label{lem-PE}
Let $E$ be a vector bundle on $X$ of rank $r$.  Denote by $\pi:\PP E
\rightarrow X$ the associated projective bundle.  The graded pieces of
the Chern character of $T_{\PP E}$ are, $c_1(T_{\PP E}) =
r\zeta + \pi^*(c_1(T_X) + c_1(E))$ and $\text{ch}_2(T_{\PP E}) =
r\zeta^2/2 + \pi^*c_1(E)\zeta +
\pi^*(\text{ch}_2(T_X)+\text{ch}_2(E))$, where $\zeta$ equals $c_1(\OO_{\PP
  E}(1))$.    
\end{lem}

\begin{proof}
There is an Euler sequence,
$$
\begin{CD}
0 @>>> \OO_{\PP E} @>>> \pi^* E\otimes \OO_{\PP E}(1) @>>> T_{\PP E/X}
@>>> 0.
\end{CD}
$$
Therefore $\text{ch}(T_{\PP E/X}) = \pi^*\text{ch}(E)e^\zeta - 1$,
i.e., 
$$
\begin{array}{c}
(r+\pi^* c_1(E) + \pi^*\text{ch}_2(E)+\dots)(1 + \zeta + \zeta^2/2 +
\dots) -1 = \\
\ \ \\
\lt[ r-1 \rt] + [ r\zeta + \pi^* c_1(E) ] + [ r\zeta^2/2 +
\pi^*c_1(E)\zeta + \pi^* \text{ch}_2(E)] + \dots
\end{array}
$$
Using the exact sequence,
$$
\begin{CD}
0 @>>> T_{\PP E/X} @>>> T_{\PP E} @>>> \pi^* T_X @>>> 0,
\end{CD}
$$
$\text{ch}(T_{\PP E})$ equals $\text{ch}(T_{\PP E/X}) +
\pi^*\text{ch}(T_X)$.  Thus $\text{ch}_1(T_{\PP E/X}) = r\zeta +
\pi^*(c_1(T_X) + c_1(E))$ and,
$$
\text{ch}_2(T_{\PP E}) = r\zeta^2/2 + \pi^*c_1(E) \zeta +
\pi^*(\text{ch}_2(T_X) + \text{ch}_2(E)).
$$
\end{proof}

\begin{prop} \label{prop-FanoPE}
Let $X$ be a smooth Fano manifold and let $E$ be a vector bundle on
$X$ of rank $r$.  The projective bundle $\PP E$ is Fano iff there
exists $\epsilon > 0$ such that for every irreducible curve $B\subset
X$, 
$$
\mu_B^1(E|_B) - \mu_B(E|_B) \leq (1-\epsilon)\text{deg}_B(-K_X)/r,
$$
where $\mu_B$ and $\mu_B^1$ are the slopes from
Definition~\ref{defn-mu}, resp. Definition~\ref{defn-mu2}.
\end{prop}

\begin{proof}
The invertible sheaf $\omega_{\PP E}^\vee$ is $\pi$-relatively ample.
By hypothesis, $\omega_X^\vee$ is ample.  By Lemma~\ref{lem-relample},
$\omega_{\PP E}^\vee$ is ample iff there exists a real number
$\epsilon >0$ such that 
$$
\text{deg}_B(g^*\omega_{\PP E}^\vee) \geq \epsilon
\text{deg}_B(g^*\pi^*\omega_X^\vee),
$$
for every finite morphism $g:B\rightarrow \PP E$ of a smooth, connected
curve to $X$ for which $\pi\circ g$ is also finite.
Using the universal property of $\PP E$, this holds iff for every
finite morphism $f:B\rightarrow X$ and every invertible quotient $f^*
E^\vee \rightarrow L^\vee$, 
$$
\text{deg}_B(g^*\omega_{\PP E}^\vee) \geq \epsilon \text{deg}_B(g^*
\pi^*\omega_X^\vee),
$$
where $g:B\rightarrow \PP E$ is the associated morphism.  
By Lemma~\ref{lem-PE}, $\text{deg}_B(\omega_{\PP E}^\vee)$ equals
$rc_1(L^\vee) + c_1(f^*E) + c_1(f^* T_X)$, i.e.,
$$
r[c_1(f^* T_X)/r - (\mu_B(L)-\mu_B(f^*E))].
$$
So, finally, $\omega_{\PP E}^\vee$ is ample iff there exists $\epsilon
>0$ such that for every finite morphism $f:B\rightarrow X$ and every
invertible quotient $f^*E^\vee \twoheadrightarrow L^\vee$,
$$
\mu_B(L) - \mu_B(f^* E) \leq (1-\epsilon)\text{deg}_B(f^*c_1(T_X))/r.
$$
Taking the supremum over covers of $B$ and invertible quotients of the
pullback of $E$, this
is,
$$
\mu_B^1(f^* E) - \mu_B(f^* E) \leq (1-\epsilon)\text{deg}_B(-f^*K_X)/r.
$$
Since every finite morphism $f:B\rightarrow X$ factors through its image,
it suffices to consider only irreducible curves $B$ in $X$. 
\end{proof}

\medskip\noindent
For $r=2$, there is a
necessary and sufficient condition for $\text{ch}_2(T_{\PP E})$ to be nef.

\begin{prop} \label{prop-rk2}
Let $E$ be a vector bundle on $X$ of rank $2$.  Denoting by $\pi:\PP E
\rightarrow X$ the projection, $\text{ch}_2(T_{\PP E}) = \pi^*
(\text{ch}_2(T_X)+ 1/2(c_1^2-4c_2)(E) )$.  Therefore
$\text{ch}_2(T_{\PP E})$ is nef iff
$\text{ch}_2(T_X)+1/2(c_1^2-4c_2)(E)$ 
is nef.  If $\text{dim}(X)>0$, $\text{ch}_2(T_{\PP
  E})$ is not weakly positive.
\end{prop}

\begin{proof}
By Lemma~\ref{lem-PE}, $\text{ch}_2(T_{\PP E})$ equals $\zeta^2 +
\pi^*c_1(E) \zeta + \pi^*(\text{ch}_2(T_X)+\text{ch}_2(E))$.  By
definition of the Chern classes of $E$, $\zeta^2 + \pi^*c_1(E) \zeta +
\pi^* c_2(E)$ equals $0$.  So the class above is $-\pi^*c_2(E) +
\pi^*(\text{ch}_2(T_X) + \text{ch}_2(E))$.  Finally, $\text{ch}_2(E)
-c_2(E)$ equals $1/2(c_1^2 -2c_2)(E) - c_2(E) = 1/2(c_1^2-4c_2)(E)$.  
\end{proof}

\medskip\noindent
Applying Proposition~\ref{prop-FanoPE} and Proposition~\ref{prop-rk2}
to the vector bundle $E=L^\vee \oplus \OO_X$ gives Example 2 in
Section~\ref{sec-nefexamples}.  

\medskip\noindent
Finally, for $r>2$, 
there is a necessary condition for $\text{ch}_2(T_{\PP E})$ to be nef.

\begin{prop} \label{prop-ch2P}
Let $E$ be a vector bundle of rank $r>2$ on $X$.  If $\text{ch}_2(T_{\PP E})$
is nef, then the pullback of $E$ to every smooth, projective,
connected
curve is semistable.  Also,
$\text{ch}_2(T_{\PP E})$ is not weakly positive if
$\text{dim}(X) > 0$ and if the
pullback of $E$ to some curve is strictly semistable, e.g., if $X$
contains a rational curve. 
\end{prop}

\begin{proof}
If the pullback of $E$ to some smooth, projective, connected curve is
not semistable, then by Corollary ~\ref{cor-mu2}, there exists a
smooth, projective, connected curve $B$, a morphism $f:B\rightarrow
X$, and a rank 2 locally free subsheaf $F$ of $f^*E$ such that
$f^*E/F$ is locally free and $\mu_B(F) > \mu_B(E)$.  There is an
induced morphism $g:\PP F \rightarrow \PP E$ such that $\pi\circ g =
f\circ \pi$.  By Lemma~\ref{lem-PE}, $g^*\text{ch}_2(T_{\PP E})$
equals $r\xi^2/2 + \pi^*f^*c_1(E)\xi +
\pi^*f^*(\text{ch}_2(T_X)+\text{ch}_2(E))$, where $\xi$ equals
$c_1(\OO_{\PP F}(1))$.  Since $B$ is a curve, $f^*(\text{ch}_2(T_X) +
\text{ch}_2(E))$ equals $0$.  Also, by definition of the Chern classes
of $F$, $\xi^2 + \pi^*c_1(F) \xi = 0$.  Substituting in,
$$
g^* \text{ch}_2(T_{\PP E}) = 1/2\pi^*(2c_1(f^*E) -rc_1(F))\xi.
$$
In particular, $\text{deg}_{\PP F}(g^* \text{ch}_2(T_{\PP E}))$ equals
$1/2(2\text{deg}_B(c_1(f^*E)) - r\text{deg}_B(F))$.  This equals 
$r(\mu_B(f^*E) - \mu_B(F))$, which is negative by construction.
Therefore $\text{ch}_2(T_{\PP E})$ is not nef.  
\end{proof}

\begin{rmk} \label{rmk-ch2P}
A vector bundle on a product of projective spaces whose restriction to
every curve is semistable is of the form $L^{\oplus r}$, where $L$ is
an invertible sheaf, ~\cite[Thm. 3.2.1]{OSS}.  In this case, $\PP E$
is also a product of projective spaces.  
\end{rmk}

\begin{cor} \label{cor-ch2P}
Let $X$ be a Fano manifold.  For every vector bundle $E$ on $X$ of
rank $r>1$, $\text{ch}_2(T_{\PP E})$ is not weakly positive.
\end{cor}

\section{Blowings up} \label{sec-blowup}

\medskip\noindent
Let $X$ be a smooth, connected, projective variety, let $i:Y\hookrightarrow X$ 
be the closed immersion of a smooth, connected  
subvariety of $X$ of codimension $c$.  
Denote by $\nu:\widetilde{X} \rightarrow X$ the
blowing up of $X$ along $Y$.  Denote by $\pi:E\rightarrow Y$ the
exceptional divisor.  Denote by $j:E\rightarrow \widetilde{X}$ the
obvious inclusion.
Then $E = \PP N_{Y/X}$ and
$i^*\OO_{\widetilde{X}}(E)$ is canonically isomorphic to $\OO_{\PP
  N}(-1)$.  

\begin{lem} \label{lem-bup}
The graded pieces of the Chern character of $\widetilde{X}$ are,
$c_1(T_{\widetilde{X}}) = \nu^* c_1(T_X) - (c-1)[E]$ and
$\text{ch}_2(T_{\widetilde{X}}) = \nu^*\text{ch}_2(T_X) + (c+1)[E]^2/2 -
i_* \pi^* c_1(N_{Y/X})$
\end{lem}

\begin{proof}
Using the short exact sequence,
$$
\begin{CD}
0 @>>> \nu^* \Omega_X @>>> \Omega_{\widetilde{X}} @>>> j_*\Omega_\pi
@>>> 0,
\end{CD}
$$
$\text{ch}(\Omega_{\widetilde{X}})$ equals $\nu^*\text{ch}(\Omega_X) +
\text{ch}(j_*\Omega_\pi)$.  Grothendieck-Riemann-Roch for the morphism
$j$ gives,
$$
\text{ch}(Rj_* a) = j_*(\text{ch}(a))(1-e^{-[E]})/[E].
$$
Using the Euler sequence for $\Omega_\pi$,
$$
\begin{CD}
0 @>>> \Omega_\pi @>>> \pi^* N_{Y/X}^\vee \otimes \OO_{\PP N}(-1) @>>>
\OO_E @>>> 0,
\end{CD}
$$
$\text{ch}(\Omega_\pi)$ equals
$\pi^*\text{ch}(N_{Y/X}^\vee)i^*(1+e^{[E]}) -1.
$
Putting the pieces together gives the lemma.
\end{proof}

\medskip\noindent
When is $\widetilde{X}$ Fano?  Denote by $\mc{C}_1$ the collection of
finite morphisms $g:B\rightarrow X$ from a smooth, connected curve to
$X$ whose image is not contained in $Y$.  Denote by $\mc{C}_2$ the
collection of finite morphisms $g:B\rightarrow Y$ from a smooth,
connected curve to $Y$.  The following result is well-known.

\begin{prop} \label{prop-bup}
Let $h$ be the first Chern class of an ample invertible sheaf on $X$,
e.g., $h = c_1(T_X)$ if $X$ is Fano.  The blowing up $\widetilde{X}$
is Fano iff there exists $\epsilon > 0$ such that, 
\begin{enumerate}
\item[(i)]
for every $g:B\rightarrow X$ in $\mc{C}_1$, 
$$
\text{deg}_B(g^{-1}Y)
\leq\frac{1}{c-1} (\text{deg}_B(g^* c_1(T_X)) - \epsilon
\text{deg}_B(g^*h)),
$$
and
\item[(ii)]
for every $g:B\rightarrow Y$ in $\mc{C}_2$, 
$$
\mu_B^1(g^*N_{Y/X}) \leq
\frac{1}{c-1}(\text{deg}_B(g^*c_1(T_X)) - \epsilon
\text{deg}_B(g^*h)).
$$
\end{enumerate}
\end{prop}

\medskip\noindent
The proof is similar to the proof of Proposition ~\ref{prop-FanoPE}.
Using an analogue of Proposition ~\ref{prop-rk2}, no blowing-up of
$\PP^n$ is a Fano manifold with $\text{ch}_2$ nef.  

\section{Theorems about vector bundles on curves}

\medskip\noindent
There are two theorems in this section.  The first
goes back to Shou-Wu Zhang, though possibly
it is older.  
The second is a variation of the first.  

\begin{defn} \label{defn-const}
Let $B$ be a smooth, projective curve.  A \emph{cover} of $B$ is a
finite, flat morphism $f:C\rightarrow B$ of constant, positive
degree.  A \emph{vector bundle} on $B$ is a locally free
$\OO_B$-module of constant rank.
\end{defn}

\begin{defn} \label{defn-mu}
Let $B$ be a smooth,  
projective curve.  For every non-zero vector bundle
$E$ on $B$, the \emph{slope} is,
$$
\mu_B(E) = \text{deg}(E)/\text{rank}(E) = \chi(B,E)/\text{rank}(E) -
\chi(B,\OO_B). 
$$
For every cover
$f:C\rightarrow B$ 
and every non-zero vector bundle
$E$ on $C$, the \emph{$B$-slope} is,
$$
\mu_{B}(f,E) :=
\text{deg}(E)/(\text{deg}(f)\text{rank}(E)) = \mu_B(f_*E)-\mu_B(f_*\OO_C).
$$
When there is no chance of confusion, this is denoted simply $\mu_B(E)$.
\end{defn}

\medskip\noindent
For every cover $g:C'\rightarrow C$, $f\circ g:C'\rightarrow B$ is a
cover and 
$\mu_B(f\circ g, g^*E)$ equals $\mu_B(f,E)$.  

\begin{defn} \label{defn-mu2}
Let $B$ be a smooth, projective curve and let $E$ be a
vector bundle on $B$ of rank $r>0$.  For every integer $1\leq k\leq
r$, define $\mu^k_B(E)$ to be,
$$
\begin{array}{c}
\sup \{ -\mu_B(f,F^\vee) | f:C\rightarrow B \text{ a cover },
f^*E^\vee \rightarrow 
F^\vee \text{ a rank } k \text{ quotient} \} \\
\ \ \\
= \sup\{ \mu_B(f,F) | f:C\rightarrow B \text{ a cover }, F\subset f^*E
\text{ a rank } k \\
\text{ subbundle whose cokernel is locally free} \}.
\end{array}
$$
\end{defn}

\medskip\noindent
Let $f:X\rightarrow Y$ be a morphism of projective varieties.  Denote
by $\mc{C}_1$ the collection of all irreducible curves in $X$ not
contained in a fiber of $f$.  Denote by $\mc{C}_2$ the collection of
finite morphisms $g:C\rightarrow X$ occurring as the normalization of
an irreducible curve in $X$ not contained in a fiber of $f$.  Finally,
denote by $\mc{C}_3$ the collection of all finite morphisms from
smooth, connected curves to $X$ whose image is not contained in a
fiber of $f$.

\begin{lem} \label{lem-relample}
Let $f:X\rightarrow Y$ be a morphism of projective varieties and let $L$
be an ample invertible $\OO_Y$-module.  
An $f$-ample
invertible $\OO_X$-module $M$ is ample iff 
there exists a real number $\epsilon > 0$ such that
for every morphism
$g:C\rightarrow X$ in $\mc{C}_1$, resp. $\mc{C}_2,\mc{C}_3$, 
$\text{deg}_C(g^*M) \geq \epsilon
\text{deg}_C(g^*f^*L)$.  
\end{lem}

\begin{proof}
Because $M$ is $f$-ample and $L$ is ample, there exists an integer
$n>0$ such that $M\otimes f^*L^{\otimes n}$ is ample.  By Kleiman's
criterion, $M$ is ample iff there exists a real number $0< \delta <1$
such that for every irreducible curve $C$ in $X$,
$$
\text{deg}_C(M) \geq \delta \text{deg}_C(M\otimes f^*L^{\otimes n}).
$$
Simplifying, this is equivalent to,
$$
\text{deg}_C(M) \geq \frac{n\delta}{1-\delta} \text{deg}_C(f^* L).
$$
As $M$ is $f$-ample, this holds if
$C$ is contained in a fiber of $f$.  So $M$ is ample iff the
inequality holds for every curve in $\mc{C}_1$.
Setting $\epsilon =
n\delta/(1-\delta)$, $\delta = \epsilon/(n+\epsilon)$, gives the lemma.

\medskip\noindent
Since $\mc{C}_2\subset \mc{C}_3$, the condition for $\mc{C}_3$ implies
the condition for $\mc{C}_2$.  Since degrees on a curve can be
computed after pulling back to the normalization, the condition for
$\mc{C}_2$ implies the condition for $\mc{C}_1$.  Finally,
for every morphism $g:C\rightarrow X$
in $\mc{C}_3$, $g(C)$ is in $\mc{C}_1$.
The inequality for $g(C)$ implies the inequality
for $C$.  Thus the condition for $\mc{C}_1$ implies the condition for
$\mc{C}_3$. 
\end{proof}

\begin{lem} \label{lem-epsilon}
Let $B$ be a smooth, connected, projective curve.
A nonzero vector bundle $E$ on $B$ is ample iff there exists a positive
real number $\delta$ such that for every cover $f:C\rightarrow B$ and
every invertible quotient $f^*E\rightarrow L$, $\mu_B(L) \geq
\delta$.  In other words, $E$ is ample iff $\mu^1_B(L^\vee) < 0$.  
\end{lem}

\begin{proof}
Denote by $\pi:\PP E^\vee \rightarrow B$ the projective bundle
associated to $E^\vee$, and denote by $\pi^*E \rightarrow \OO_{\PP
  E^\vee}(1)$ the tautological invertible quotient.  By definition,
$E$ is ample iff $\OO_{\PP E^\vee}(1)$ is an ample invertible sheaf.
Of course $\OO_{\PP E^\vee}(1)$ is $\pi$-relatively ample.  Let $M$ be
an invertible $\OO_B$-module of degree $1$.  Then $M$ is ample.  
By Lemma~\ref{lem-relample},
$\OO_{\PP E^\vee}(1)$ is ample iff there exists $\epsilon > 0$ such
that for every smooth, connected curve $C$ and every finite morphism
$g:C\rightarrow \PP E^\vee$ such that $\pi\circ g$ is finite, 
$\text{deg}_C(g^*\OO_{\PP E^\vee}(1)) \geq \epsilon
\text{deg}_C(g^*\pi^*M)$.  Of course
$\text{deg}_C(g^*\pi^*M)=\text{deg}(\pi\circ g)$.  Using the universal
property of $\PP E^\vee$, this holds iff for every cover
$f:C\rightarrow B$ and every invertible quotient $f^*E\rightarrow L$,
$$
\text{deg}_C(L) \geq \epsilon \text{deg}(f) \Leftrightarrow \mu_B(L)
\geq \epsilon.
$$
\end{proof}

\begin{lem} \label{lem-Zhang}
For every ample vector bundle $E$ on $B$, there exists a cover
$f:C\rightarrow B$, invertible $\OO_C$-modules $L_1,\dots,L_r$, and a
morphism of $\OO_C$-modules, $\phi:f^*E\rightarrow (L_1\oplus \dots
\oplus L_r)$ such that,
\begin{enumerate}
\item[(i)]
the support of $\text{coker}(\phi)$ is a finite set,
\item[(ii)]
for every $i=1,\dots,r$, the projection $f^*E \rightarrow
\oplus_{j\neq i} L_j$ is surjective, and
\item[(iii)]
for every $i=1,\dots,r$, $\mu_B(L_i) = \text{deg}_B(E)$.
\end{enumerate}
\end{lem}

\begin{proof}
Denote $r=\text{rank}(E)$.  The claim is that for every $k=1,\dots,r$,
there exists a cover $f_k:C_k\rightarrow B$, invertible
$\OO_{C_k}$-modules $L_{k,1},\dots,L_{k,k}$, and a morphism of
$\OO_{C_k}$-modules, $\phi_k:f^*E \rightarrow (L_{k,1}\oplus \dots \oplus
L_{k,k})$ satisfying (ii) and (iii) above and the following variant of
(i): for $k<r$, $\phi_k$ is surjective and for $k=r$, the support of
$\text{coker}(\phi_k)$ is a finite set.
The lemma is the case
$k=r$.  The claim is proved by induction on $k$.

\medskip\noindent
The base case is $k=1$.  Denote by $\pi:\PP E^\vee \rightarrow B$ the
projective bundle associated to $E^\vee$, and denote by
$\pi^* E \twoheadrightarrow \OO_{\PP E^\vee}(1)$ the tautological
invertible quotient.  By hypothesis, $\OO_{\PP E^\vee}(1)$ is ample.
By Bertini's theorem, for $d_1,\dots,d_{r-1} \gg 0$, there exist
effective Cartier divisors $D_1,\dots,D_{r-1}$ with $D_i\in |\OO_{\PP
  E^\vee}(d_i)|$ such that the intersection $C_1 = D_1\cap \dots \cap
D_r$ is a smooth, connected curve, cf. ~\cite{Jou}.  Denote by
$f_1:C_1\rightarrow B$ the restriction of $\pi$.  Denote by
$\phi_1:f^*E \rightarrow L_{1,1}$ the restriction of $\pi^* E \rightarrow
\OO_{\PP E^\vee}(1)$.  This satisfies (i) because $\pi^* E \rightarrow
\OO_{\PP E^\vee}(1)$ is surjective.  It satisfies (ii) trivially.
Finally, $\text{deg}(f)$ equals $d_1\times \cdots
\times d_{r-1}$, and $\text{deg}_{C_1}(L_{1,1})$ equals $d_1\times \cdots
\times d_{r-1} \times [c_1(\OO_{\PP E^\vee}(1))]^r$, i.e., $d_1\times
\cdots \times d_{r-1}\times \text{deg}_B(E)$.  Therefore $\mu_B(L_{1,1}) =
\text{deg}_B(E)$, i.e., this satisfies (iii).

\medskip\noindent
By way of induction, assume the result is known for $k<r$, and
consider the case $k+1$.  Since $\phi_k$ is surjective, there is an
induced closed immersion $\PP(L_{k,1}\oplus\dots\oplus L_{k,k})^\vee
\hookrightarrow \PP (f_k^* E)^\vee$.  The image is irreducible and has
codimension $r-k\geq 1$.  For every $i=1,\dots,k$, the image of
$\PP(\oplus_{j\neq i} L_{k,j})^\vee$ is irreducible and has codimension
$r-k+1\geq 2$.  Associated to the finite morphism $f_k$, there is a
finite morphism $\PP (f_k^* E)^\vee \rightarrow \PP E^\vee$.  The
pullback of an ample invertible sheaf by a finite morphism is ample;
hence $\OO_{\PP (f_k^* E)^\vee}(1)$ is ample.  By Bertini's theorem,
for $d_1,\dots,d_{r-1} \gg 0$, there exist effective Cartier divisors
$D_1,\dots, D_{r-1}$ with $D_i \in |\OO_{\PP (f_k^* E)^\vee}(d_i)|$
such that the intersection $C_{k+1} = D_1\cap \dots \cap D_{r-1}$ is a
smooth, connected curve, disjoint from $\PP(\oplus_{j\neq i}
L_j)^\vee$ for every $i=1,\dots,k$, and either disjoint from
$\PP(\oplus_i L_i)^\vee$ if $k<r-1$, or else intersecting
$\PP(\oplus_i L_i)^\vee$ in finitely many points if $k=r-1$.
Define $g_{k+1}:C_{k+1}
\rightarrow C_k$ to be the restriction of the projection.  Define
$f_{k+1} = f_k \circ g_{k+1}$, define $L_{k+1,i} = g_{k+1}^* L_{k,i}$
for $i=1,\dots,k$, and define $L_{k+1,k+1}$ to be the restriction of
$\OO_{\PP (f_k^* E)^\vee}(1)$.  Define $\phi_{k+1}$ to be the obvious
morphism. 

\medskip\noindent
The cokernel of $\phi_{k+1}$ is supported on the intersection of
$C_{k+1}$ with $\PP(L_{k,1}\oplus \dots \oplus L_{k,k})^\vee$.  By
construction, this is empty if $k<r-1$, and is a 
a finite set if $k=r-1$.  Thus $\phi_{k+1}$ satisfies (i).
By the induction hypothesis, $f_{k+1}^* E\rightarrow
(L_{k+1,1}\oplus\dots \oplus L_{k+1,k})$, which is the pullback under
$g_{k+1}$ of $\phi_k$, is surjective.  For $i=1,\dots,k$, the cokernel
of $f_{k+1}^* E\rightarrow \oplus_{j\neq i} L_{k+1,j}$ is supported on
the intersection of $C_{k+1}$ with the image of $\PP(\oplus_{j\neq i}
L_{k,j})^\vee)$.  By construction, this is empty, i.e., $f_{k+1}^* E
\rightarrow \oplus_{j\neq i} L_{k+1,j}$ is surjective.  Thus
$\phi_{k+1}$ satisfies (ii).  Finally, $\phi_{k+1}$ satisfies (iii) by
the same argument as in the base case.  The claim is proved by
induction on $k$.
\end{proof}

\begin{thm} \label{thm-epsilon}
For every non-zero vector bundle $E$ on $B$, for
every $\epsilon > 0$, there exists a cover
$f:C\rightarrow B$ and a invertible quotient
$f^*E \rightarrow L$ such that $\mu_B(L) < \mu_B(E) + \epsilon$.
In other words, $\mu_B^1(E^\vee) \geq \mu_B(E^\vee)$.
\end{thm}

\begin{proof}
Denote $r=\text{rank}(E)$.
If $r=1$, set $f=\text{Id}_B$ and $L=E$.  Then $L$ is an invertible
quotient of $f^* E$, and 
$\mu_B(L)$ equals $\mu_B(E)$ which is less than $\mu_B(E)+\epsilon$.
Therefore assume $r > 1$.

\medskip\noindent
Certainly an effective version of the following argument can be given,
but a simpler argument is by contradiction.

\begin{hyp} \label{hyp-cont}
For every cover
$f:C\rightarrow B$ and every invertible quotient $f^*E\rightarrow L$,
$\mu_B(L)$ is $\geq \mu_B(E) + \epsilon$, i.e., $\mu_B^1(E^\vee) <
\mu_B(E^\vee) - \epsilon$.    
\end{hyp}

\medskip\noindent
By way of contradiction, assume Hypothesis~\ref{hyp-cont}.
Let $f:C\rightarrow
B$ be a connected, smooth cover of degree $d$.
For every $a/d \in \frac{1}{d}\ZZ$, there exists an
invertible sheaf $M$ on $C$ of degree $a$, and thus $\mu_B(M) = a/d$.
In particular, for $d$ sufficiently large, there exists an invertible
quotient $M$ such that $0 < \mu_B(E) - \mu_B(M) < \epsilon/(r-1)$.
Denote $\delta = \mu_B(E) -\mu_B(M)$.   
Denote $F = f^*E\otimes M^\vee$.  Then
$\mu_B(F)$ equals $\delta$, and $0<\delta < \epsilon/(r-1)$.  

\medskip\noindent
Let $g:C'\rightarrow C$ be any cover and let $g^*F \rightarrow N$ be
any invertible quotient.  Then $f\circ g:C'\rightarrow B$ is a cover
and $(f\circ g)^* E = g^*F\otimes g^*M \rightarrow N\otimes g^*M$ is
an invertible quotient.  By Hypothesis~\ref{hyp-cont},
$$
\begin{array}{c}
\mu_C(N) = \text{deg}(f)\mu_B(N) =\text{deg}(f)(\mu_B(N\otimes g^*M) -
\mu_B(M)) \\
\\
\geq \text{deg}(f)((\mu_B(E) + \epsilon) -\mu_B(M)) >
\text{deg}(f)\epsilon.
\end{array}
$$
By Lemma~\ref{lem-epsilon}, $F$ is an ample vector bundle on $C$.  By
Lemma~\ref{lem-Zhang}, there exists a cover $g:C'\rightarrow C$ and an
invertible quotient $g^*F\rightarrow P$ such that $\mu_B(P) =
r\mu_B(F) = r\delta$.
Therefore $L:= g^*M\otimes
P$ is an invertible quotient of $g^*f^* E$
and,
$$
\mu_B(L) = \mu_B(g^*M\otimes P) = \mu_B(M) +
r\delta = 
\mu_B(E) + (r-1)\delta.
$$
By hypothesis, $(r-1)\delta < \epsilon$.  So $\mu_B(L) < \mu_B(E) +
\epsilon$, contradicting Hypothesis~\ref{hyp-cont}.  
The proposition is proved by contradiction.  
\end{proof}

\begin{cor} \label{cor-epsilon}
For every non-zero vector bundle $E$ on $B$, for every $\epsilon > 0$,
there exists a cover $f:C\rightarrow B$ and a sequence of vector
bundle quotients,
$$
f^*E = E^r \twoheadrightarrow E^{r-1} \twoheadrightarrow \dots
\twoheadrightarrow E^1,
$$
such that each $E^k$ is a vector bundle of rank $k$ and $\mu_B(E^k) <
\mu_B(E) + \epsilon$.
\end{cor}

\begin{proof}
The proof is by induction on the rank $r$ of $E$.  If $\text{rank}(E)=1$,
defining $f=\text{Id}_B$ and $E^1 = E$, the result follows.  
Thus, assume
$r>1$ and the result is known for smaller values of $r$.
By Theorem~\ref{thm-epsilon}, there exists a cover
$g:B'\rightarrow B$ and a rank 1 quotient $g^*E\rightarrow L$
such that $\mu_B(L) < \mu_B(E) + \epsilon$.  Denote by $K$ the
kernel of $g^*E\rightarrow L$.  Then
$\text{rank}(K)=r-1$ and $\mu_B(K) = (r\mu_B(E) - \mu_B(L))/(r-1)$.
By the induction hypothesis, there
exists a cover $h:C\rightarrow B'$ and a sequence of vector bundle
quotients,
$$
h^* K = K^{r-1} \twoheadrightarrow \dots \twoheadrightarrow K^1,
$$
such that each $K^k$ is a vector bundle of rank $k$, and
$\mu_{B'}(K^k) \leq \mu_{B'}(K) + \text{deg}(g)\epsilon$.  
Of course $\mu_B(F) =
\mu_{B'}(F)/\text{deg}(g)$ for every $F$.  Thus $\mu_B(K^k) \leq
\mu_B(K) + \epsilon$.  

\medskip\noindent
Define $f=h\circ g$, define
$E^1 = h^* L$, and for every $k=2,\dots, r$, define $f^* E
\twoheadrightarrow E^k$ to be the unique quotient whose kernel
is contained in $h^*K$ and such that $h^*K\rightarrow E^k$
has image $K^{k-1}$.  Then $\mu_B(E^1) = \mu_B(L) \leq
\mu_B(E)+\epsilon$, and for $k=2,\dots, r$, 
$$
\begin{array}{c}
\mu_B(E^k) = 1/k(\mu_B(L) + (k-1)\mu_B(K^{k-1})) < 1/k(\mu_B(L) +
(k-1)\mu_B(K) + (k-1)\epsilon) = \\
\\
\frac{r(k-1)}{(r-1)k}\mu_B(E) + \frac{r-k}{(r-1)k}\mu_B(L) +
\frac{(r-1)(k-1)}{(r-1)k}\epsilon < \mu_B(E) +
\frac{r-k}{(r-1)k}\epsilon + \frac{(r-1)(k-1)}{(r-1)k}\epsilon <
\mu_B(E) + \epsilon. 
\end{array}
$$
\end{proof}

\medskip\noindent
For semistable bundles in characteristic zero, 
there is a more precise result.

\begin{thm}[Zhang] \label{thm-Zhang}
Let $B$ be a smooth, projective curve over an algebraically closed
field of characteristic $0$.  Let $E$ be a semistable vector bundle on
$B$.  
Let
$\epsilon$ be a positive real number.
There exists a cover $f:C\rightarrow B$, invertible sheaves
$L_1,\dots,L_r$ on $C$, and a morphism of $\OO_C$-modules,
$\phi:f^*E\rightarrow (L_1\oplus \dots \oplus L_r)$ such that,
\begin{enumerate}
\item[(i)]
the support of 
$\text{coker}(\phi)$ is a finite set,
\item[(ii)]
for every $i=1,\dots,r$, the projection $f^*E \rightarrow
\oplus_{j\neq i} L_j$ is surjective,
\item[(iii)]
for every $i=1,\dots,r$, $\mu_B(L_i) \leq \mu_B(E)+\epsilon$.
\end{enumerate}
\end{thm}

\begin{proof}
Denote $r=\text{rank}(E)$.  If $r$ equals $1$, the theorem is
trivial.  Thus assume $r>1$.  
As in the proof of Theorem~\ref{thm-epsilon}, there exists a cover
$g:C'\rightarrow B$ and an invertible sheaf $M$ on $C'$ such that
$0<\mu_B(E) - \mu_B(M) < \epsilon/(r-1)$.  Denote $\delta = \mu_B(E)
-\mu_B(M)$ and denote $F=g^*E\otimes M^\vee$.  Then $\mu_B(F)$ equals
$\delta$, and $0<\delta < \epsilon/(r-1)$.  

\medskip\noindent
Let $h:C\rightarrow C'$ be any cover and let $h^*F\rightarrow N$ be an
invertible quotient.  The composition $g\circ h:C\rightarrow B$ is a
cover.  By Kempf's theorem, ~\cite{Kempf}, which ultimately relies on
the theorem that every stable vector bundle admits a Hermite-Einstein
metric, $(g\circ h)^* E$ is semistable.  (Note, there are
counterexamples in positive characteristic.)  Therefore $h^* F$ is
semistable.  So $\mu_{C}(L) \geq \mu_{C}(h^* F)$, i.e., $\mu_{C'}(L)
\geq \mu_{C'}(F) =\delta$.  Thus by Lemma~\ref{lem-epsilon}, $F$ is an ample
vector bundle on $C'$.  Thus by Lemma~\ref{lem-Zhang}, there exists a
cover $h:C\rightarrow C'$, invertible $\OO_{C}$-modules
$N_1,\dots,N_r$, and a morphism of $\OO_{C}$-modules $\psi:h^*F
\rightarrow (N_1\oplus \dots \oplus N_r)$ satisfying (i), (ii) and
(iii) of Lemma~\ref{lem-Zhang}.  Define $f = g\circ h$, $L_i = N_i
\otimes h^* M$ and $\phi$ is the twist of $\psi$ by
$\text{Id}_{h^*M}$.  Then $\phi$ satisfies (i) and (ii).  And for
every $i=1,\dots,r$, 
$$
\begin{CD}
\mu_B(L_i) = \mu_B(N_i) + \mu_B(M) =
\mu_{C'}(N_i)/\text{deg}(g) + \mu_B(E) - \delta = \\
\mu_B(E) +
r\delta/\text{deg}(g) - \delta \leq \mu_B(E) +
(r-1)\delta/\text{deg}(g) < \mu_B(E) + \epsilon.
\end{CD}
$$
\end{proof}

\medskip\noindent
Of course, $\mu^r_B(E)$ equals $\mu_B(E)$.  
The other values are more interesting.

\begin{cor} \label{cor-mu2}
The slopes $\mu_B^k(E)$ satisfy $\mu_B^1(E)\geq \mu_B^2(E) \geq \dots
\geq \mu_B^r(E) = \mu_B(E)$.  For each $1\leq k < r$,
$\mu_B^k(E)=\mu_B(E)$ 
iff $f^*E$ is semistable for every cover $f:C\rightarrow B$.
\end{cor}

\begin{proof}
By Corollary~\ref{cor-epsilon}, for every $\epsilon >0$, there exists
a cover $f:C\rightarrow B$ and a rank $k$ quotient $f^*E\rightarrow
E^k$ such that $\mu_B(E^k) < \mu_B(E)+\epsilon$.  Thus $\mu^k_B(E)\geq
\mu_B(E)$.  Applying the same reasoning to rank $k-1$ quotients of
rank $k$ quotients of $f^*E$, $\mu^{k-1}_B(E)\geq \mu^k_B(E)$.  

\medskip\noindent
If $f^* E$ is semistable for every cover $f:C\rightarrow B$, then
every vector bundle quotient of $f^* E$ has slope $\geq \mu_C(f^*E)$,
and thus has $B$-slope $\geq \mu_B(f^*E)$.  Therefore $\mu^k_B(E) \leq
\mu_B(E)$, i.e., $\mu^k_B(E) = \mu_B(E)$.  

\medskip\noindent
Conversely, suppose there is a cover $f:C\rightarrow B$ such that $f^*
E$ is not semistable.  Then there exists a vector bundle quotient
$f^*E \twoheadrightarrow F$ such that $\mu_B(F)<\mu_B(E)$.  Denote the rank
by $l$.  Suppose first that $l\geq k$, and define $\epsilon =
\text{deg}(f)(\mu_B(E)-\mu_B(F))$.  Then by Corollary
~\ref{cor-epsilon}, there exists a cover $g:C'\rightarrow C$ and a
rank $k$ quotient $g^*F \twoheadrightarrow G$ such that $\mu_C(G) < \mu_C(F)
+ \epsilon$.  Therefore $g^*f^*E \twoheadrightarrow g^*F
\twoheadrightarrow G$ is a rank $k$ quotient of $g^*f^*E$ and
$\mu_B(G) < \mu_C(F) + (\mu_B(E)-\mu_B(F)) = \mu_B(E)$.  Therefore
$\mu^k_B(E) > \mu_B(E)$.  

\medskip\noindent
Next suppose that $l<k$.  Denote by $K$ the kernel of $f^*E\rightarrow
F$.  Then $r\mu_B(E) = l\mu_B(F) + (r-l)\mu_B(K)$.  Define,
$$
\epsilon = \frac{(r-k)l\text{deg}(f)(\mu_B(E)-\mu_B(F))}{(r-l)(k-l)}.
$$
By Corollary ~\ref{cor-epsilon}, there exists a cover $g:C'\rightarrow
C$ and a rank $k-l$ quotient $g^*K \twoheadrightarrow G'$ such that
$\mu_C(G') < \mu_C(K) + \epsilon$.  Therefore $\mu_B(G') < \mu_B(K) +
\epsilon/\text{deg}(f)$.  Define $g^*f^*E\rightarrow G$ to be the
unique vector bundle whose kernel is contained in $g^*K$ and such that
the image of $g^*K\rightarrow G$ equals $G'$.  Then,
$$
\begin{array}{c}
k\mu_B(G) = l\mu_B(F) + (k-l)\mu_B(G') < l\mu_B(F) + (k-l)\mu_B(K) +
(k-l)\epsilon/\text{deg}(f) = \\
\\
l\mu_B(F) + \frac{k-l}{r-l}(r\mu_B(E) - l\mu_B(F)) +
\frac{k-l}{\text{deg}(f)}\epsilon = \\
\\
k\mu_B(E) -
\frac{(r-k)l}{r-l}(\mu_B(E) - \mu_B(F)) +
\frac{(r-k)l}{r-l}(\mu_B(E)-\mu_B(F)) = k\mu_B(E).
\end{array}
$$
Thus $\mu_B(G) < \mu_B(E)$, and therefore $\mu^k_B(E)>\mu_B(E)$.  
\end{proof}

\bibliography{my}
\bibliographystyle{alpha}

\end{document}